\documentclass{article}

\usepackage{epsfig,xspace}
\usepackage{latexsym}
\usepackage[matrix,arrow,curve]{xy}\xyoption{all}

\usepackage[usenames,dvipsnames,svgnames,table]{xcolor}
\usepackage{amssymb,amsmath,latexsym,amsfonts,amsthm,alltt}
\usepackage[usenames,dvipsnames,svgnames,table]{xcolor}
\usepackage[dvips,usenames]{stmaryrd}
\usepackage{enumitem}
\usepackage{url}
\usepackage{graphicx}
\usepackage{float}
\usepackage{amstext}
\usepackage{cite}
\usepackage{hyperref}
\usepackage{MnSymbol}

\newtheorem{defn}{Definition}[section]
\newtheorem{thm}[defn]{Theorem}
\newtheorem{lemma}[defn]{Lemma}
\newtheorem{prop}[defn]{Proposition}
\newtheorem{coro}[defn]{Corollary}

\newcommand{\ra}{\rightarrow}
\newcommand{\lra}{\longrightarrow}

\newcommand{\Ra}{\Rightarrow}
\newcommand{\La}{\Leftarrow}

\newcommand{\llra}{\longleftrightarrow}

\newcommand{\midsp}{\;|\;}
\newcommand{\sub}[2]{#1_{{}_{#2}}}
\newcommand{\telos}{\hfill$\Box$}
\newcommand{\type}[1]{{\tt #1}}
\newcommand{\iso}{\backsimeq}

\newcommand{\val}[1]{\mbox{$\lsem #1\rsem$}}
\newcommand{\forces}{\Vdash}
\newcommand{\dforces}{\forces^{\!\!\partial}}
\newcommand{\yvval}[1]{\mbox{$\llparenthesis #1 \rrparenthesis $}}

\newcommand{\upv}{\upVdash}
\newcommand{\rperp}{\mbox{${}^{\upv}$}}
\newcommand{\gphi}{{\mathcal  G}(Y)}
\newcommand{\gpsi}{{\mathcal  G}(X)}

\newcommand{\lperp}{{}\rperp}

\DeclareMathOperator{\filt}{Filt}
\DeclareMathOperator{\idl}{Idl}

\newcommand{\Mtright}{\mathrel{\mbox{$|\!\!\largetriangleright$}}}
\newcommand{\Mtleft}{\mathrel{\mbox{$\largetriangleleft\!\!|$}}}

\title{Choice-free Topological Duality\\ for Implicative Lattices
 and Heyting Algebras
}
\author{Chrysafis (Takis) Hartonas\\
University of Thessaly, Greece\\
$\type{hartonas@uth.gr}$}

\begin{document}
\maketitle

\begin{abstract}
We develop a common semantic framework for the interpretation both of $\mathbf{IPC}$, the intuitionistic propositional calculus, and of logics weaker than $\mathbf{IPC}$ (substructural and subintuitionistic logics). This is done by proving a choice-free representation and duality theorem for implicative lattices, which may or may not be distributive. The duality specializes to a choice-free duality for the category of Heyting algebras and a category of topological sorted frames with a ternary sorted relation.
\end{abstract}

\section{Introduction}
Stone \cite{stone2}, Priestley \cite{hilary,hilary_davey_2002} and Esakia \cite{Esakia2019-ESAHA} duality for distributive lattices and Heyting algebras, in particular, have served as the basic background for the extension of relational semantics to the case of distributive logic systems. Esakia duality more specifically, establishes an equivalence between Heyting and $\mathbf{S4}$-algebras (closure algebras) \cite[Section~3.4]{Esakia2019-ESAHA}.  Esakia's \cite{esakia} and Blok's \cite{blok} independently discovered results on the isomorphism of the lattice of superintuitionistic logics and the lattice of normal extensions of $\mathbf{S4}$ led to the development of the Blok-Esakia theory, brief reviews of which can be found in \cite{chagrov,Wolter2014}. 

The present article focuses  not on superintuitionistic (intermediate) logics, but instead on logics weaker than $\mathbf{IPC}$, the intuitionistic propositional calculus, including both substructural logics \cite{ono-galatos}, as well as subintuitionistic logics \cite{Visser1981,dosen-modal-trans-in-KD,corsi-subint,restall-subint,Wansing1997-subint,celani-subint-2001,weak-subint}. Building on previous work by the author on duality for normal lattice expansions \cite{duality2} and on a recently advanced research program of choice-free topological dualities, initiated by Bezhanishvili and Holliday \cite{choice-free-BA} and soon followed by similar results \cite{choice-free-Ortho, choice-free-dmitrieva-bezanishvili,choice-free-deVries,choice-free-inquisitive-DBLP:conf/wollic/BezhanishviliGH19,choiceFreeStLog}, we present in this article a sorted relational framework for the interpretation of the language of $\mathbf{IPC}$, in which both subintuitionistic and the implicational fragment of substructural logics can be interpreted. 

We prove a choice-free topological duality for implicative lattices and, in particular, for Heyting algebras. The duality specializes to the intermediate case of distributive implicative lattices. The distinctive feature of any subintermediate logic (substructural, or subintuitionistic) is that implication is no longer residuated with conjunction and then distribution may be either forced into the axiomatization (as in the case of subintuitionistic logics), or left out (as is the common case for substructural logics). We also show by a representation argument that every implicative lattice is a reduct of a residuated lattice, and that every distributive implicative lattice is a reduct of a residuated Heyting algebra, where these are structures similar to the J\'{o}nsson and Tsinakis residuated Boolean algebras \cite{residBA}, but with a Heyting rather than a Boolean underlying algebra.

In devising a common semantic (relational) framework for the interpretation of $\mathbf{IPC}$ and of weaker logics, the semantics of intuitionistic implication remains the standard one
\begin{equation}\label{int clause1}
w\forces\varphi\ra\psi\;\mbox{ iff }\; \forall u(wRu\;\mbox{ $\lra$}\;(u\forces\varphi\lra u\forces \psi))
\end{equation}
but for weaker logics a ternary relation is to be used (and the two clauses are equivalent in the case of $\mathbf{IPC}$). This is a setting familiar from the Routley-Meyer semantics \cite{routley1973semantics} for Relevance Logic, with relevant implication interpreted by a clause of the form
\begin{equation}\label{relevant clause}
w\forces\varphi\ra\psi\;\mbox{ iff }\; \forall u,z(u\forces\varphi\mbox{ $\wedge$ }wRuz\;\mbox{ $\lra$ }\; z\forces\psi).
\end{equation}
An information-theoretic interpretation has been proposed in the context of Barwise's channel theory \cite{Barwise1993-BARCCA,lif}, where implication sentences designate properties of information channels $w$, connecting sites $u,z$ (written $u\stackrel{w}{\leadsto}z$) and the clause takes the form $w\forces\varphi\ra\psi\;\mbox{ iff }\; \forall u,z(u\forces\varphi\mbox{ and }u\stackrel{w}{\leadsto}z\;\mbox{ implies }\; z\forces\psi)$.

Dropping distribution, as well, relational semantics for the resulting logic systems has been given in sorted frames (polarities) $\mathfrak{P}=(X,\upv,Y)$ with additional relations \cite{mai-gen,mai-grishin,COUMANS201450,suzuki8,discres,pnsds,redm}, referred to as information systems, object-attribute systems or formal contexts in the formal concept analysis (FCA) tradition \cite{wille2}, or polarities, after Birkhoff \cite{birkhoff}. The relation ${\upv}\subseteq X\times Y$ induces a Galois connection $(\;)\rperp:\powerset(X)\leftrightarrows\powerset(Y):{}\rperp(\;)$ in the standard way and sentences are interpreted as Galois stable sets $A={}\rperp(A\rperp)$. A co-interpretation (refutation) set is also associated to a sentence and thus the semantics uses both a relation of satisfaction, $\forces$, and one of co-satisfaction (refutation), $\dforces$, and the clause for implication may be written in the form
\begin{equation}\label{implication sat clause}
x\forces\varphi\ra\psi \;\mbox{ iff }\; \forall u,z(u\forces\varphi\mbox{ and }xRuz\;\mbox{ implies }\; z\not\dforces\psi)
\end{equation}
We show that, in the case of $\mathbf{IPC}$, implication can be equivalently interpreted by a clause of the form in \eqref{int clause1}.

When distribution is not assumed and the logic (lattice) does not come equipped with a De Morgan negation (complementation) operator, interpreting disjunction becomes an additional issue needing attention. In \cite{choice-free-dmitrieva-bezanishvili}, the semantics of disjunction is based on a type of lattice completion that the authors call an $F^2$-completion, which is an iterated filter construction, similar to the one used by Gouveia and Priestley \cite{hilary-sem} to investigate the problem of the canonical extension of a semilattice. The resulting clause is along the lines of \eqref{F2disj}
\begin{equation}\label{F2disj}
x\forces\varphi\vee\psi\;\mbox{ iff }\; (x\forces\varphi\mbox{ or }x\forces\psi), \;\mbox{ or }\;
\exists u,z(u\forces\varphi\mbox{ and }z\forces\psi\mbox{ and }u\wedge z\leq x)
\end{equation}
where the carrier set of the frame is assumed to have a semilattice structure. We proceed differently, resorting to what was coined `order-dual relational semantics' in \cite{odigpl,odjpl} (where, however, general frames only were used), dually interpreting disjunction as its order-dual in a sorted frame $(X,\upv, Y)$, namely as conjunction, by clause \eqref{disj clause}, where $x\in X$ and $y\in Y$,
\begin{equation}\label{disj clause}
x\forces\varphi\wedge\psi\mbox{ iff }x\forces\varphi\mbox{ and }x\forces\psi\hskip1cm y\dforces\varphi\vee\psi\mbox{ iff }y\dforces\varphi\mbox{ and }y\dforces\psi
\end{equation}
an intuitive reading of which is that $y$ refutes a disjunction $\varphi\vee\psi$ iff it refutes both disjuncts $\varphi$ and $\psi$. Interpretation $\val{\vartheta}$ and cointerpretation $\yvval{\vartheta}$ of a sentence $\vartheta$ are related by $\val{\vartheta}={}\rperp\yvval{\vartheta}$ and $\yvval{\vartheta}=\val{\vartheta}\rperp$ hence, if desired, a satisfaction clause for disjunction can be stated: $x\forces\varphi\vee\psi\;\mbox{ iff }\; \forall y(xIy\;\lra\; (y\dforces\varphi\;\lra\;y\not\dforces\psi))$, where $I$ is the complement of the relation $\upv$ of the sorted frame (polarity).
\\[1mm]

In this article, though applications in logic motivate what we do, we do not discuss logic matters directly and we lay, instead, the groundwork for a common semantic framework for $\mathbf{IPC}$ and for weaker logics in the language of $\mathbf{IPC}$.
\\[1mm]

Section \ref{impl lattice section} is a preparatory section where we merely define implicative lattices and residuated Heyting algebras, as particular cases of interest of normal lattice expansions, so that we can relate to results obtained in \cite{duality2,choiceFreeStLog}. 

Section \ref{impl frames section} is devoted to a study of sorted frames $\mathfrak{F}=(X,\upv,Y,T)$, with a ternary relation $T\subseteq Y\times(X\times Y)$. Frame basics are reviewed in Section \ref{frame prelim section}.  Implicative frames are defined (axiomatized) and their full complex algebra of Galois stable sets is studied in Section \ref{complex algebra section}. Section \ref{distr section} establishes a first-order condition for the full complex algebra of the frame to be (completely) distributive. Finally, Section \ref{Heyting frames section} defines Heyting frames, by suitable first-order conditions, whose full complex algebra is a complete Heyting algebra and, moreover, the implication operator on stable sets induced by the ternary relation $T$ is identical to the residual of intersection.

Section \ref{rep section} is devoted to a choice-free representation and duality for implicative lattices, specializing it to distributive implicative lattices and Heyting algebras. In Section \ref{rep integral section} a representation of integral implicative lattices is detailed, including the case of an underlying distributive lattice. Section \ref{rep Heyting section} turns to the representation of Heyting algebras. Duality is discussed in Section \ref{duality section}, which relies heavily on both \cite{duality2,choiceFreeStLog}. 

Pointers to applications and  concluding remarks are given in Section \ref{conc section}.

\section{Implicative Lattices}
\label{impl lattice section}
Let $\{1,\partial\}$ be a 2-element set, $\mathbf{L}^1=\mathbf{L}$ and $\mathbf{L}^\partial=\mathbf{L}^{\mathrm{op}}$ (the opposite lattice). Extending established terminology \cite{jt1}, a function $f:\mathbf{L}_1\times\cdots\times\mathbf{L}_n\lra\mathbf{L}_{n+1}$ will be called {\em additive} and {\em normal}, or a {\em normal operator}, if it distributes over finite joins of the lattice $\mathbf{L}_i$, for each $i=1,\ldots n$, delivering a join in $\mathbf{L}_{n+1}$.

\begin{defn}\rm
\label{normal lattice op defn}
An $n$-ary operation $f$ on a bounded lattice $\mathbf{L}$ is {\em a normal lattice operator of distribution type  $\delta(f)=(i_1,\ldots,i_n;i_{n+1})\in\{1,\partial\}^{n+1}$}  if it is a normal additive function  $f:\mathbf{L}^{i_1}\times\cdots\times\mathbf{L}^{i_n}\lra\mathbf{L}^{i_{n+1}}$ (distributing over finite joins in each argument place), where  each $i_j$, for  $j=1,\ldots,n+1$,   is in the set $\{1,\partial\}$, hence $\mathbf{L}^{i_j}$ is either $\mathbf{L}$, or $\mathbf{L}^\partial$.

If $\tau$ is a tuple (sequence) of distribution types, a {\em normal lattice expansion of (similarity) type $\tau$} is a lattice with a normal lattice operator of distribution type $\delta$ for each $\delta$ in $\tau$.

The {\em category {\bf NLE}$_\tau$}, for a fixed similarity type $\tau$, has normal lattice expansions of type $\tau$ as objects. Its morphisms are the usual algebraic homomorphisms.
\end{defn}

In this article we focus on bounded lattices $\mathbf{L}=(L,\leq,\wedge,\vee,0,1,\ra)$ with a binary operation $\ra$ of implication, of distribution type $\delta(\ra)=(1,\partial;\partial)$, defined below. As $\tau$ consists of a single distribution type, we let $\mathbf{NLE}_{1\partial\partial}$ designate the category of integral implicative lattices. The objects of $\mathbf{NLE}_{1\partial\partial}$ are defined below and its morphisms are the usual algebraic bound-preserving homomorphisms.

\begin{defn}\rm\label{implicative lattice defn}
An {\em integral implicative lattice}   $\mathbf{L}=(L,\leq,\wedge,\vee,0,1,\ra)$ is a bounded lattice with a binary (implication) operation where the following axioms hold, in addition to the axioms for bounded lattices.
\begin{tabbing}
\hskip5mm\=(A1)\hskip5mm\= $(a\vee b)\ra c=(a\ra c)\wedge(b\ra c)$\\
\>(A2)\> $a\ra(b\wedge c)=(a\ra b)\wedge(a\ra c)$\\
\>(A3)\> $a\leq b\;\llra \; 1\leq a\ra b$
\end{tabbing}
For example, the implication-only fragment of an $\mathbf{FL}_{ew}$-algebra (a residuated lattice assuming exchange $a\circ b=b\circ a$ and weakening $a\circ b\leq b$) is an integral implicative lattice.

The lattice is  {\em distributive} if axiom (A4) is assumed.
\begin{tabbing}
\hskip5mm\=(A4)\hskip5mm\= $a\wedge(b\vee c)=(a\wedge b)\vee(a\wedge c)$
\end{tabbing}
It is a Heyting algebra if axioms (H1)-(H2) are added to the axioms for bounded lattices
\begin{tabbing}
\hskip5mm\=(H1)\hskip5mm\= $a\wedge(a\ra b)\leq b$\\
\>(H2)\> $b\leq a\ra(a\wedge b)$
\end{tabbing}
\end{defn}
By $\mathbb{IIL}$ we designate the quasi-variety of integral implicative lattices. Similarly, d$\mathbb{IIL}\subset\mathbb{IIL}$ designates the quasi-variety of distributive integral implicative lattices. Heyting algebras are distributive integral implicative lattices, hence the variety $\mathbb{HA}$ of Heyting algebras is contained in the quasi-variety d$\mathbb{IIL}$. 

\begin{lemma}\rm
Countably many distinct quasi-varieties $\mathbb{V}$ are contained in the complete lattice of quasi-varieties between $\mathbb{HA}$ and $\mathbb{IIL}$.
\end{lemma}
\begin{proof}
Define $a^n\ra b$ recursively by  $a^0\ra b=b$ and $a^{n+1}\ra b=a\ra(a^n\ra~b~)$. For $n\geq 1$, let $A_n$ be the axiom $a^{n+1}\ra b\leq a^n\ra b$ and let $\mathbb{V}_n$ be the subvariety of $\mathbb{IIL}$ generated by adding the axiom $A_n$ to $\mathbb{IIL}$, so that $\mathbb{V}_n\subset\mathbb{IIL}$, for all $n\geq 1$. Notice that $\mathbb{V}_1$ adds the contraction axiom $a^2\ra b=a\ra(a\ra b)\leq a\ra b$ to $\mathbb{IIL}$, which is intuitionistically valid, hence $\mathbb{HA}\subset\mathbb{V}_1$. In $\mathbb{V}_n=\mathbb{IIL}+A_n$ we have $a^{n+2}\ra b=a\ra(a^{n+1}\ra b)\leq a\ra(a^n\ra b)=a^{n+1}\ra b$, so $A_{n+1}$ is derivable from the quasi-equational theory of $\mathbb{V}_n$. Thus $\mathbb{V}_n\subseteq\mathbb{V}_{n+1}$. By contrast, $A_n$ is not derivable in $\mathbb{V}_{n+1}$, hence $\mathbb{V}_n\subset\mathbb{V}_{n+1}$. There is therefore a countable chain of varieties between $\mathbb{HA}$ and $\mathbb{IIL}$, such that $\mathbb{HA}\subset\mathbb{V}_1\subset\cdots\subset\mathbb{V}_n\subset\cdots\subset\mathbb{IIL}$.
\end{proof}

\begin{defn}\rm
\label{resid Heyting defn}
A {\em residuated Heyting algebra} is a structure \[\mathbf{H}=(H,\leq,\wedge,\vee,0,1,\ra,/,\circ,\backslash)\] where $(H,\leq,\wedge,\vee,0,1,\ra)$ is a Heyting algebra, $(H,\leq,\wedge,\vee,0,1,/,\circ,\backslash)$ is a distributive integral residuated lattice and $(H,\leq,\wedge,\vee,0,1,\backslash)$ is an integral implicative lattice.
\end{defn}
J\'{o}nsson and Tsinakis \cite{residBA} have studied residuated Boolean algebras in the context of their study of relation algebras. Residuated Heyting algebras arise in the context of representation and canonical extension of distributive integral implicative lattices (cf Proposition \ref{reduct prop}).

The main fact about implicative lattices is stated in the following theorem.
\begin{thm}\rm
\label{reduct thm}
Every integral implicative lattice $\mathbf{L}$ is a reduct of an integral residuated lattice $\mathbf{R}$. If $\mathbf{D}$ is a distributive integral implicative lattice, then it is a reduct of a residuated Heyting algebra $\mathbf{H}$. In both cases, $\mathbf{R,H}$ are constructed as a canonical extension of $\mathbf{L}$ and of $\mathbf{D}$, respectively, and they are therefore unique, up to an isomorphism that fixes $\mathbf{L}$ (respectively, $\mathbf{D}$).
\end{thm}
The proof will be provided in Section \ref{rep section}, by a representation argument (see Corollary \ref{reduct coro} and Corollary \ref{Heyting rep coro} for the special case of Heyting algebras). 

\section{Sorted Residuated Implicative Frames}
\label{impl frames section}
\subsection{Frame Preliminaries}
\label{frame prelim section}
Consider $\{1,\partial\}$ as a set of sorts and let $Z=(Z_1,Z_\partial)$ be a sorted set.
Sorted residuated frames $\mathfrak{F}=(Z_1,\upv,Z_\partial)$ are triples consisting of nonempty sets $Z_1=X,Z_\partial=Y$ and a binary relation ${\upv}\subseteq X\times Y$.
We refer to the relation $\upv$ as the {\em Galois relation} of the frame. It generates a Galois connection in the standard way $(\;)\rperp:\powerset(X)\leftrightarrows\powerset(Y)^\partial:\lperp(\;)$ ($V\subseteq U\rperp$ iff $U\subseteq\lperp V$)
\begin{tabbing}
\hskip2cm\=$U\rperp$\hskip2mm\==\hskip1mm\=$\{y\in Y\midsp\forall x\in U\; x\upv y\}$ \hskip1mm\==\hskip1mm\= $\{y\in Y\midsp U\upv y\}$\\
\>$\lperp V$\>=\>$\{x\in X\midsp \forall y\in V\;x\upv y\}$\>=\>$\{x\in X\midsp x\upv V\}$.
\end{tabbing}

A subset $A\subseteq X$ will be called {\em stable} if $A={}\rperp(A\rperp)$. Similarly, a subset $B\subseteq Y$ will be called {\em co-stable} if $B=({}\rperp B)\rperp$. Stable and co-stable sets will be referred to as {\em Galois sets}, disambiguating to {\em Galois stable} or {\em Galois co-stable} when needed and as appropriate.

By $\gpsi,\gphi$ we designate the complete lattices of stable and co-stable sets, respectively. Note that the Galois connection restricts to a dual isomorphism $(\;)\rperp:\gpsi\iso\gphi^\partial:{}\rperp(\;)$.

Preorder relations are induced on each of the sorts, by setting for $x,z\in X$, $x\preceq z$ iff $\{x\}\rperp\subseteq\{z\}\rperp$ and, similarly, for $y,v\in Y$, $y\preceq v$ iff ${}\rperp\{y\}\subseteq{}\rperp\{v\}$.
A (sorted) frame is called {\em separated} if the preorders $\preceq$ (on $X$ and on $Y$) are in fact partial orders $\leq$. All frames are hereafter assumed to be separated.

Our notational conventions are these of \cite[Remark~3.2]{duality2}. Vectorial notation $\vec{u}$ is used for a tuple of elements $(u_1,\ldots,u_n)$, for some $n$. 
We use $\Gamma$ to designate  upper closure  $\Gamma U=\{z\in X\midsp\exists x\in U\;x\leq z\}$, for $U\subseteq X$, and similarly for $U\subseteq Y$. The set $U$ is {\em increasing} (an upset) iff $U=\Gamma U$. For a singleton set $\{x\}\subseteq X$ we write $\Gamma x$, rather than $\Gamma(\{x\})$ and similarly for $\{y\}\subseteq Y$.

We typically use the standard Formal Concept Aanalysis  priming notation for each of the two Galois maps ${}\rperp(\;),(\;)\rperp$. This allows for stating and proving results for each of $\gpsi,\gphi$ without either repeating definitions and proofs, or making constant appeals to duality. Thus for a Galois set $G$, $G'=G\rperp$, if $G\in\gpsi$ ($G$ is a Galois stable set), and otherwise $G'={}\rperp G$, if $G\in\gphi$ ($G$ is a Galois co-stable set).

For an element $u$ in either $X$ or $Y$ and a subset $W$, respectively of $Y$ or $X$, we write $u|W$, under a well-sorting assumption, to stand for either $u\upv W$ (which stands for $u\upv w$, for all $w\in W$), or $W\upv u$ (which stands for $w\upv u$, for all $w\in W$), where well-sorting means that either $u\in X, W\subseteq Y$, or $W\subseteq X$ and $u\in Y$, respectively. Similarly for the notation $u|v$, where $u,v$ are elements of different sort.

Let $\mathfrak{F}=(X,\upv,Y)$ be a  polarity and $u$ a point in $Z=X\cup Y$. By \cite[Lemma~3.3]{duality2}, the following are some basic facts.
\begin{enumerate}
\item $\upv$ is increasing in each argument place (and thereby its complement $I$ is decreasing in each argument place).
\item $(\Gamma u)'=\{u\}'$ and $\Gamma u=\{u\}^{\prime\prime}$ is a Galois set.
\item Galois sets are increasing, i.e. $u\in G$ implies $\Gamma u\subseteq G$.
\item For a Galois set $G$, $G=\bigcup_{u\in G}\Gamma u$.
\item For a Galois set $G$, $G=\bigvee_{u\in G}\Gamma u=\bigcap_{v|G}\{v\}'$.
\item For a Galois set $G$ and any set $W$, $W^{\prime\prime}\subseteq G$ iff $W\subseteq G$.
\end{enumerate}

We refer to principal upper sets $\Gamma x\in\gpsi (x\in X= \filt(\mathbf{L}))$, as {\em closed}, or {\em filter} elements of $\gpsi$ and to sets ${}\rperp\{y\}\in\gpsi\; (y\in Y=\idl(\mathbf{L}))$ as {\em open}, or {\em ideal} elements of $\gpsi$, and similarly for sets $\Gamma y, \{x\}\rperp$ with $x\in X, y\in Y$. This creates an unfortunate clash with topological terminology and we shall have to rely on context to disambiguate. Furthermore, a closed element $\Gamma u$ is said to be {\em clopen} if $\Gamma u=\{w\}'$ for some $w$, which is unique, by frame separation. If $\Gamma u$ is clopen, then we call the point $u$ a {\em clopen point}, as well.

By \cite[Lemma~3.3]{duality2}, the closed elements of $\gpsi$  join-generate $\gpsi$, while the open elements meet-generate $\gpsi$ (similarly for $\gphi$).

For a sorted relation $R\subseteq\prod_{j=1}^{j=n+1}Z_{i_j}$, where $i_j\in\{1,\partial\}$ for each $j$ (and thus $Z_{i_j}=X$ if $i_j=1$ and $Z_{i_j}=Y$ when $i_j=\partial$), we make the convention to regard it as a relation $R\subseteq Z_{i_{n+1}}\times\prod_{j=1}^{j=n}Z_{i_j}$, we agree to write its sort type as $\sigma(R)=(i_{n+1};i_1\cdots i_n)$ and for a tuple of points of suitable sort we write $uRu_1\cdots u_n$ for $(u,u_1,\ldots,u_n)\in R$. If $\sigma$ is the sort of $R$ we often display it by writing $R^\sigma$. It is understood that if $\sigma_1\neq\sigma_2$, then $R^{\sigma_1}$ and $R^{\sigma_2}$ name different relations (the sort superscript is part of the name designation).

\begin{defn}[Galois dual relation] \label{Galois dual relations}\rm
For a relation $R$, of sort type $\sigma$, its {\em Galois dual relation} $R'$ is the relation defined by $uR'\vec{v}$ iff $\forall w\;(wR\vec{v}\lra w|u)$. In other words, $R'\vec{v}=(R\vec{v})'$.
\end{defn}

\begin{defn}[Sections of relations]\rm
\label{sections defn}
For an $(n+1)$-ary relation $R^\sigma$ (of sort $\sigma$) and an $n$-tuple $\vec{u}$, $R^\sigma\vec{u}=\{w\midsp wR^\sigma\vec{u}\}$ is the {\em section} of $R^\sigma$ determined by $\vec{u}$. To designate a section of the relation at the $k$-th argument place we let $\vec{u}[\_]_k$ be the tuple with a hole at the $k$-th argument place. Then $wR^\sigma\vec{u}[\_]_k=\{v\midsp wR^\sigma\vec{u}[v]_k\}\subseteq Z_{i_k}$ is the $k$-th section of $R^\sigma$.
\end{defn}

The frames of interest in the present article are the implicative frames studied in the next section.

\subsection{Implicative Frames and their Full Complex Algebras}
\label{complex algebra section}
\begin{defn}\rm
\label{implicative frame defn}
If $\mathfrak{F}=(X,\upv,Y,T^{\partial 1\partial})$ is a frame with a ternary relation $T$ of the indicated sort type, i.e. $T\subseteq Y\times (X\times Y)$, then $\mathfrak{F}$ is an {\em implicative frame} iff the axioms of Table \ref{axioms table1} hold.
\end{defn}

\begin{table}[!htbp]
\caption{Implicative Frame Axioms}
\label{axioms table1}
\begin{enumerate}
  \item[(F0)] $\forall x,y(x\upv y\llra\forall u\; uT'xy)$
  \item[(F1)] The frame is separated
  \item[(F2)] For all $x\in X$ and $v\in Y$, the section $[\_]Txv$ of the relation $T$ is a closed element of $\gphi$
  \item[(F3)] For any $y\in Y$, the binary relation $yT$ is decreasing in both argument places
  \item[(F4)] For any $z\in X$, both sections $zT'[\_]v$ and $zT'x[\_]$ of the Galois dual relation $T'$ of $T$ are Galois sets
\end{enumerate}
\end{table}

A sorted image operator $\alpha_T=\largetriangleright$ is defined on $U\subseteq X, V\subseteq Y$ by
\[
U{\largetriangleright} V=\alpha_T(U,V)=\{y\in Y\midsp\exists x,v(x\in U\;\wedge\;v\in V\;\wedge\;yTxv)\}=\bigcup_{x\in U}^{v\in V}Txv
\]
and we let ${\Mtright}:\gpsi\times\gphi\lra\gphi$ designate the closure $\overline{\alpha}_T$ of the restriction of $\alpha_T=\largetriangleright$ to Galois sets, defined on $A\in\gpsi, B\in\gphi$ by
\[
A\Mtright B=\overline{\alpha}_T(A,B)=\left(\{y\in Y\midsp\exists x,v(x\in A\;\wedge\;v\in B\;\wedge\;yTxv)\}\right)''=\bigvee_{x\in A, v\in B}Txv
\]
\begin{lemma}\rm
\label{triangle dist lemma}
The Galois set operation ${\Mtright}:\gpsi\times\gphi\lra\gphi$ distributes over arbitrary joins of Galois sets in each argument place, i.e. for $A_i\in \gpsi$ with $i\in I$ and $B_j\in\gphi$ with $j\in J$, $\left(\bigvee_{i\in I}A_i\right)\Mtright\left(\bigvee_{j\in J}B_j\right)=\bigvee_{i\in I,j\in J}(A_i\Mtright B_j)$.
\end{lemma}
\begin{proof}
  The claim is a special instance, for a ternary relation $T$, of \cite[Theorem~3.12]{duality2}, given the axioms for implicative frames in Table \ref{axioms table1}.
\end{proof}

\begin{defn}\rm
\label{implication def}
For $A,C\in\gpsi$, define $A\Ra C=(A\Mtright C')'={}\rperp(A\Mtright C\rperp)$.
\end{defn}
Notice that $A\Ra C=(A\Mtright C')'=(A\largetriangleright C')'''=(A\largetriangleright C')'$.

\begin{prop}\rm
\label{Ra prop}
The following hold
\begin{enumerate}
  \item $\left(\bigvee_{i\in I}A_i\right)\Ra\left(\bigcap_{j\in J}C_j\right)=\bigcap_{i\in I,j\in J}(A_i\Ra C_j)$
  \item $A\Ra C=\bigcap_{x\in A, C\upv y}(\Gamma x\Ra{}\rperp\{y\})$
  \item $uT'xy$ iff $u\in(\Gamma x\Ra{}\rperp\{y\})$, for all $u,x\in X$ and $y\in Y$
  \item $u\in (A\Ra C)$ iff $\forall x\in X\forall y\in Y(x\in A\;\wedge\;C\upv y\lra uT'xy)\}$
  \item $A\subseteq C$ iff $X\subseteq A\Ra C$, for any $A,C\in\gpsi$.
\end{enumerate}
\end{prop}
\begin{proof}
For (1), since the Galois connection is a duality of stable and co-stable sets, $(\;)\rperp:\gpsi\iso\gphi^\partial:{}\rperp(\;)$, it follows from Lemma \ref{triangle dist lemma}, which establishes that $\Mtright$ distributes over arbitrary joins in both argument places, that the distribution type of $\Ra$ is $(1,\partial;\partial)$, in other words
\[
\left(\bigvee_{i\in I}A_i\right)\Ra\left(\bigcap_{j\in J}C_j\right)=\bigcap_{i\in I,j\in J}(A_i\Ra C_j).
\]
For (2), closed elements join-generate and open elements meet-generate $\gpsi$, hence we obtain in particular that
\begin{equation}
A\Ra C=\left(\bigvee_{x\in A}\Gamma x\right)\Ra\left(\bigcap_{C\upv y}{}\rperp\{y\}\right)=\bigcap_{x\in A, C\upv y}(\Gamma x\Ra{}\rperp\{y\})\nonumber
\end{equation}
For claim (3), by definition $\Gamma x\Ra{}\rperp\{y\}={}\rperp(\Gamma x\Mtright\Gamma y)=\left(\bigvee_{x\leq z,y\leq v}Tzv\right)'=\bigcap_{x\leq z,y\leq v}T'zv$. By stability of the sections of $T'$ and since stable sets are increasing, the set $T'xy$ is contained in every set $T'zv$ for $x\leq z$ and $y\leq v$. It follows that $\Gamma x\Ra{}\rperp\{y\}=T'xy$ and this proves claim (3). Thereby, $A\Ra C=\bigcap_{x\in A,C\upv y}T'xy$, from which claim (4) follows.

For claim (5), suppose first that $A\subseteq C$. By claim (4) $u\in A\Ra C$ iff for every $x\in A$ and every $y\in Y$ such that $C\upv y$, we have $uT^{11\partial}xy$. Now by $x\in A\subseteq C\upv y$, it follows that $x\upv y$ and then by (F0) we get that for every $u\in X$ we have $uT^{11\partial}xy$, hence for every $u\in X$ we obtain $u\in A\Ra C$, which is to say that $X\subseteq A\Ra C$.

Conversely, if $X\subseteq A\Ra C=\{u\in X\midsp \forall x,y(x\in A\;\wedge\; C\upv y\;\lra\; uT^{11\partial}xy)\}$, then $\forall u\in X\forall x,y(x\in A\;\wedge\; C\upv y\;\lra\; uT^{11\partial}xy)$. Pushing the quantifier on $u$ we get the equivalent form $\forall x,y(x\in A\;\wedge\;C\upv y\;\lra\;\forall u\in X\;uT^{11\partial}xy)$. By (F0) this is further equivalent to $\forall x,y(x\in A\;\wedge\;C\upv y\;\lra x\upv y)$. Equivalently, we have $\forall y\in Y(C\upv y\lra A\upv y)$ which means that for any $y$ if $C\subseteq{}\rperp\{y\}$, then $A\subseteq{}\rperp\{y\}$.

Thereby
$\{y\in Y\midsp C\subseteq{}\rperp\{y\}\}\subseteq\{y\in Y\midsp A\subseteq{}\rperp\{y\}\}$, hence using meet-density of open elements we get $A=\bigcap\{{}\rperp\{y\}\midsp A\subseteq{}\rperp\{ y\}\}\subseteq\bigcap\{{}\rperp\{y\}\midsp C\subseteq{}\rperp\{ y\}\}=C$.
\end{proof}
Consequently, we have the following result.
\begin{coro}\rm
For any implicative frame $\mathfrak{F}=(X,\upv,Y,T^{\partial 1\partial})$ as in Definition \ref{implicative frame defn}, its full complex algebra $\mathfrak{F}^+=(\gpsi,\subseteq,\bigcap,\bigvee,\emptyset,X,\Ra)$ is an integral implicative lattice.\telos
\end{coro}

We show, in addition, that the full complex algebra $\mathfrak{F}^+$ of an implicative frame is a (complete) residuated lattice. This will be used in proving Theorem \ref{reduct thm}, in Section \ref{rep section}.

\begin{defn}\rm
\label{derived relations defn}
Let $T=T^{\partial 1\partial}$ be the frame relation, of the indicated sort, i.e. $T\subseteq Y\times (X\times Y)$. Additional ternary relations, derived from $T$, are defined by
\begin{tabbing}
\hskip2mm\=$T^{11\partial}$ \hskip4mm\= Galois dual relation of $T^{\partial 1\partial}$ \hskip4mm\= $xT^{11\partial}zv$ iff $\forall y\in Y(yT^{\partial 1\partial}zv\lra x\upv y)$\\[1mm]
\>$R^{\partial 11}$ \> argument permutation \> $vR^{\partial 11}zx$ iff $xT^{11\partial}zv$\\
\>$R^{111}$ \> Galois dual relation of $R^{\partial 11}$ \> $uR^{111}zx$ iff $\forall v(vR^{\partial 11}zx\lra u\upv v)$\\[1mm]
\>$S^{\partial\partial 1}$ \> argument permutation \> $yS^{\partial\partial 1}vx$ iff $yT^{\partial 1\partial} xv$\\
\>$S^{1\partial 1}$ \> Galois dual relation of $S^{\partial\partial 1}$ \> $uS^{1\partial 1}vx$ iff $\forall y(yS^{\partial\partial 1}vx\lra u\upv y)$.
\end{tabbing}
\end{defn}

Let $\alpha_R=\bigodot$ be the sorted image operator generated by the relation $R=R^{111}$ and set $\overline{\alpha}_R=\bigovert$ to designate the closure of its restriction to Galois stable sets.

Let also $\alpha_S=\largetriangleleft$ be the sorted image operator generated by $S$ and let furthermore $\Mtleft{=} \;\overline{\alpha}_S$ be the closure of the restriction of $\alpha_S$ to Galois sets.

\begin{defn}
  \rm 
  For $A,C\in\gpsi$, define $C\La A=(C'\Mtleft A)'={}\rperp(C\rperp\La A)$.
\end{defn}

\begin{prop}\rm
\label{residuation in frame prop}
For any Galois stable sets $A,F,C\in\gpsi$ we have
\begin{enumerate}
\item $C\La F=\{z\in X\midsp\forall x\in X\forall v\in Y(x\in F\;\wedge\;C\upv v\lra xS^{1\partial 1}vz)\}$
\item $A\subseteq C\La F$ iff $A\bigovert F\subseteq C$ iff $F\subseteq A\Ra C$
\item $\left(\bigcap_{j\in J}C_j\right)\La\left(\bigvee_{i\in I}A_i\right)=\bigcap_{i\in I,j\in J}(C_j\La A_i)$
\item $C\La A=\bigcap_{x\in A, C\upv y}({}\rperp\{y\}\La\Gamma x)$
\item $xS^{1\partial 1}yz$ iff $z\in({}\rperp\{y\}\La\Gamma x)$
\end{enumerate}
\end{prop}
\begin{proof}
The first claim is established in a way completely analogous to part 4 of Proposition \ref{Ra prop}. The second claim is verified by the following computations.
\begin{tabbing}
\hskip8mm\= $F\subseteq A\Ra C$ \\
iff\> $F\subseteq \{x\in X\midsp \forall z\in X\;\forall v\in Y\;(z\in A\;\wedge\;C\upv v\;\lra\; xT^{ 11\partial }zv))\}$\\
iff\> $\forall x\in X\;(x\in F\;\lra\; \forall z\in X\;\forall v\in Y\;(z\in A\;\wedge\;C\upv v\;\lra\;xT^{ 11\partial }zv))$\\
iff\> $\forall x,z\in X\;\forall v\in Y\;(x\in F\;\lra\; (z\in A\;\wedge\;C\upv v\;\lra\;xT^{ 11\partial }zv))$\\
iff\> $\forall x,z\in X\;\forall v\in Y\;(x\in F\;\lra\; (z\in A\;\lra\;(C\upv v\;\lra\;xT^{ 11\partial }zv)))$\\
iff\> $\forall x,z\in X\;\forall v\in Y\;(x\in F\;\wedge\; z\in A\;\lra\;(C\upv v\;\lra\;xT^{ 11\partial }zv))$\\
iff\> $\forall x,z\in X\;(x\in F\;\wedge\; z\in A\;\lra\;(C\rperp\;\subseteq\;xT^{ 11\partial }z))$\\
iff\> $\forall x,z\in X\;(x\in F\;\wedge\; z\in A\;\lra\;(C\rperp\;\subseteq\;R^{\partial 11 }zx))$\\
iff\> $\forall x,z\in X\;(x\in F\;\wedge\; z\in A\;\lra\;(\lperp(R^{\partial 11 }zx)\;\subseteq\;C))$\\
iff\> $\forall x,z\in X\;(x\in F\;\wedge\; z\in A\;\lra\;(R^{111}zx\;\subseteq\; C))$\\
iff\> $\forall u,x,z\in X\;(x\in F\;\wedge\; z\in A\;\lra\;(uR^{111}zx\;\lra\; u\in C))$\\
iff\> $\forall u,x,z\in X\;(uR^{111}zx\;\wedge\;z\in A\;\wedge\;x\in F\;\lra\;u\in C)$\\
iff\> $\forall u\;(\exists x,z\in X\;(uR^{111}zx\;\wedge\;z\in A\;\wedge\;x\in F)\;\lra\;u\in C)$\\
iff\> $\forall u\;(u\in A\bigodot F\;\lra u\in C)$\\
iff\> $A\bigodot F\subseteq C$ \\
iff\> $A\bigovert F\subseteq C$.
\end{tabbing}
Notice that, given definitions, $S'vx=S^{1\partial 1}vx=(S^{\partial\partial 1}vx)'=(T^{\partial 1\partial}xv)'=T^{11\partial}xv=T'xv$. Copying from the above computation we  have that
\begin{tabbing}
\hskip8mm\= $A\bigovert F\subseteq C$\\
iff \> $\forall x,z\in X\;\forall v\in Y\;(x\in F\;\wedge\; z\in A\;\lra\;(C\upv v\;\lra\;xT^{ 11\partial }zv))$\\
iff \> $\forall x,z\in X\;\forall v\in Y\;(z\in A\lra(x\in F\;\wedge\; C\upv v\lra xS^{1\partial 1}vz))$\\
iff \> $\forall z\in X\;(z\in A\lra \forall x\in X\forall v\in Y(x\in F\;\wedge\;C\upv v\lra xS^{1\partial 1}vz)$\\
iff \> $A\subseteq\{z\in X\midsp\forall x\in X\forall v\in Y(x\in F\;\wedge\;C\upv v\lra xS^{1\partial 1}vz)\}$\\
iff \> $A\subseteq C\La F$.
\end{tabbing}

The third claim follows by the distribution properties of $\La$, given residuation. The fourth claim is a special instance, writing $A$ as the join of $\Gamma x$, with $x\in A$, and $C$ as the meet of ${}\rperp\{y\}$, with $C\upv y$.

Finally, given definitions we have that $xS^{1\partial 1}yz$ iff $xT^{11\partial}zy$, hence $S'=S^{1\partial 1}$ is increasing in every argument place, from which the last claim follows, using claim (4).
\end{proof}
Proposition \ref{residuation in frame prop} establishes the following result.
\begin{coro}\rm
\label{resid in frame coro}
If $\mathfrak{F}=(X,\upv,Y,T^{\partial 1\partial})$ is an implicative frame, then its full complex algebra $\mathfrak{F}^+=(\gpsi,\subseteq,\bigcap,\bigvee,\emptyset,X,\La,\bigovert,\Ra)$  is an integral residuated lattice.\telos
\end{coro}

\subsection{Distributive Implicative Frames}
\label{distr section}
This section contributes by establishing a first-order condition, first specified in \cite{choiceFreeStLog}, for the lattice $\gpsi$ of stable sets to be distributive. When this is the case, intersection distributes over arbitrary joins in $\gpsi$, hence it has a residual $\widetilde{\Ra}$.  Thereby $\gpsi$ is a (complete) residuated Heyting algebra (Definition \ref{resid Heyting defn}). Furthermore, first-order conditions are also established for the case where $\bigovert$ is identified with intersection (and then also $\La,\Ra$ are identified with $\widetilde{\Ra}$).

We let $R_\leq$ be the ternary upper bound relation on $X$ defined by $xR_\leq uz$ iff both $u\leq x$ and $z\leq x$. 

\begin{prop}\rm
\label{upper bound rel prop}
Let $\mathfrak{F}=(X,\upv,Y)$ be a sorted frame (a polarity) and $\gpsi$ the complete lattice of stable sets.   If all sections of the Galois dual relation $R'_\leq$ of $R_\leq$ are Galois sets, then $\gpsi$ is completely distributive.
\end{prop}
\begin{proof}
The proof was given in \cite[Proposition~3.7]{choiceFreeStLog} which, however, has not yet appeared in print, so we repeat it here.

Let $\alpha_R$ be the image operator generated by $R_\leq$, $\alpha_R(U,W)=\bigcup_{u\in U}^{w\in W}Ruw$. Notice that, for stable sets $A,C$ (more generally, for increasing sets), $\alpha_R(A,C)=A\cap C$. Hence $\overline{\alpha}_R(A,C)=\alpha_R(A,C)=A\cap C$, since Galois sets are closed under intersection. Given the section stability hypothesis for the Galois dual relation $R'_\leq$ of $R_\leq$,  Theorem~3.12 of \cite{duality2} applies, from which distribution of $\overline{\alpha}_R$ (i.e. of intersection) over arbitrary joins of stable sets is concluded.
\end{proof}
Hence we have established the following result.
\begin{coro}\rm
\label{resid Heyting coro}
Let $\mathfrak{F}=(X,\upv,Y,T^{\partial 1\partial})$ be an implicative frame. If all sections of the Galois dual relation of the upper bound relation $R_\leq$ on $X$  are Galois sets, then the full complex algebra $\mathfrak{F}^+$ of the frame is a (complete) residuated Heyting algebra.
\end{coro}
\begin{proof}
  Combine Proposition \ref{upper bound rel prop} and Corollary \ref{resid in frame coro}.
\end{proof}

\subsection{Heyting Frames}
\label{Heyting frames section}
\begin{prop}\rm
\label{heyting prop}
Let $\mathfrak{F}=(X,\upv,Y,T^{\partial 1\partial})$ be an implicative-frame, $\mathfrak{F}^+$ its full complex algebra and $R_\leq$  the upper bound relation on $X$.  Let also $R^{111}$ be the derived relation defined in Definition \ref{derived relations defn} and assume $x,z\in X$ are arbitrary. Then the following hold:
\begin{enumerate}
  \item $\Gamma x\cap\Gamma z\subseteq\Gamma x\bigovert\Gamma z$ iff $R_\leq xz\subseteq R^{111}xz$
  \item $\Gamma x\bigovert\Gamma z\subseteq\Gamma x\cap\Gamma z$ iff $R^{111}xz\subseteq R_\leq xz$
  \item $\mathfrak{F}^+$  is a complete Heyting algebra (where $\Ra$ is residuated with intersection) iff $R^{111}xz= R_\leq xz$.
\end{enumerate}
\end{prop}
\begin{proof}
For (1), observe  that for any $y\in Y$
\begin{tabbing}
$y\in(\Gamma x\bigodot\Gamma z)'$\hskip3mm\=iff\hskip2mm\= $y\in\bigcap_{x\leq x_1}^{z\leq z_1}R'x_1z_1$\hskip3cm\= ($R$ is $R^{111}$)\\
\hskip12mm\=iff\hskip3mm\=  $\forall x_1,z_1\in X(x\leq x_1\;\wedge\;z\leq z_1\lra yR'x_1z_1)$ \> ($R'=R^{\partial 11}$)\\
\>iff\> $\forall x_1,z_1\in X(x\leq x_1\;\wedge\;z\leq z_1\lra z_1T'x_1y)$ \> (definition of $R$)\\
\>iff\> $zT'xy$   \hskip0.8cm ($T'$ sections are Galois sets, which are upsets)\\
\>iff\> $yR^{\partial 11}xz$ \>(definition of $R$)\\
\>iff\> $yR'xz$ \> ($R$ is $R^{111}$)
\end{tabbing}

We then have
\begin{tabbing}
$\Gamma x\cap\Gamma z\subseteq\Gamma x\bigovert\Gamma z$ iff $\forall u(uR_\leq xz\lra u\in\Gamma x\bigovert\Gamma z)$\\
\hskip12mm\=iff\hskip3mm\= $\forall u(uR_\leq xz\lra \Gamma u\subseteq\Gamma x\bigovert\Gamma z)$ \hskip5mm\= (Gallois sets are upsets)\\
\>iff\> $\forall u(uR_\leq xz\lra \Gamma u\subseteq\bigcap\{{}\rperp\{y\}\midsp\Gamma x\bigovert\Gamma z\subseteq{}\rperp\{y\}\}$\\
\>iff\> $\forall u(uR_\leq xz\lra\forall y(\Gamma x\bigovert\Gamma z\subseteq{}\rperp\{y\}\lra u\upv y))$\\
\>iff\> $\forall u(uR_\leq xz\lra\forall y(\Gamma y\subseteq(\Gamma x\bigovert\Gamma z)'\lra u\upv y))$\\
\>iff\> $\forall u(uR_\leq xz\lra\forall y(y\in(\Gamma x\bigodot\Gamma z)'\lra u\upv y))$\\
\>iff\> $\forall u(uR_\leq xz\lra\forall y(yR^{\partial 11}xz\lra u\upv y))$\\
\>iff\> $\forall u(uR_\leq xz\lra uR^{111}xz)$\\
\>iff\> $R_\leq xz\subseteq R^{111}xz$
\end{tabbing}

For (2), note first that for any $x,z\in X$, $R_\leq xz=\Gamma x\cap\Gamma z$, so that $R_\leq xz$ is a stable set. Note also that with $R=R^{111}$ and by Definition \ref{derived relations defn}, $R^{111}xz=(R^{\partial 11}xz)'=(zT^{11\partial}x[\;])'$. We assume (by axioms F3,F4 in Table \ref{axioms table1}) that all sections of the Galois dual relation $T^{11\partial}$ of the frame ternary relation $T^{\partial 1\partial}$ are stable, hence $R'xz=(zT^{11\partial}x[\;])''=zT^{11\partial}x[\;]=R^{\partial 11}xz$ and then $R''xz=(R^{\partial 11}xz)'=R^{111}xz=Rxz$, using the definition of $R^{111}$ as the Galois dual of $R^{\partial 11}$. The proof of claim (2) is now by the following computation.
\begin{tabbing}
$\Gamma x\bigovert\Gamma z\subseteq\Gamma x\cap\Gamma z$ iff $\Gamma x\bigodot\Gamma z\subseteq R_\leq xz$\\
\hskip13mm\= iff\hskip2mm\= $R'_\leq xz\subseteq(\Gamma x\bigodot\Gamma z)'$\\
\>iff\> $\forall y(yR'_\leq xz\lra y\in(\Gamma x\bigodot\Gamma z)'$\\
\>iff\> $\forall y(yR'_\leq xz\lra yR'xz)$\hskip1cm\=($R=R^{111}$, $R'=R^{\partial 11}, R''=R$)\\
\>iff\> $\forall y(yR''xz\lra yR''_\leq xz)$\\
\>iff\> $\forall y(yR^{111}xz\lra yR_\leq xz)$ \> ($R_\leq xz$ is a stable set)\\
\>iff\> $R^{111}xz\subseteq R_\leq xz$
\end{tabbing}

For (3), if $\mathfrak{F}^+$ is a complete Heyting algebra, then $\bigovert=\cap$, by uniqueness of adjoints and Proposition \ref{residuation in frame prop}. Using (1) and (2), it follows that  $R^{111}xz= R_\leq xz$.

For the converse, assuming that for any $x,z\in X$ we have $ R^{111}xz= R_\leq xz$, it suffices to prove that $\bigovert=\cap$, which means that $\Ra$ is the residual of $\cap$, given Proposition \ref{residuation in frame prop}.

By residuation of $\bigovert$ with $\Ra$ we obtain that $A\bigovert F=\bigvee_{x\in A,z\in F}\Gamma x\bigovert\Gamma z$. Hence for any stable set $C$, $A\bigovert F\subseteq C$ iff for all $x\in A$ and $z\in F$ we have $\Gamma x\bigovert\Gamma z\subseteq C$. Since by meet-density of open elements $C=\bigcap_{C\upv y}{}\rperp\{y\}$ we obtain that $A\bigovert F\subseteq C$ iff for all $x\in A$, $z\in F$ and $y\in Y$ such that $C\upv y$ it holds that $\Gamma x\bigovert\Gamma z\subseteq{}\rperp\{y\}$. Taking $C=\Gamma x\cap\Gamma z$, notice first that $\Gamma x\cap\Gamma z\subseteq{}\rperp\{y\}$ iff $yR'_\leq xz$, where $R'_\leq$ is the Galois dual relation of the upper bound relation $R_\leq$. Furthermore, $\Gamma x\bigovert\Gamma z\subseteq{}\rperp\{y\}$ iff $y\in(\Gamma x\bigodot\Gamma z)'$ iff  $yR'xz$, where $R'=R^{\partial 11}$ is the Galois dual of the relation $R^{111}$ (and $yR^{\partial 11}xz$ iff $zT^{11\partial}xy$ holds, by definition). Thereby we have
\begin{tabbing}
$A\bigovert F\subseteq A\cap F$\hskip5mm\= iff \hskip3mm\=$\forall x\in A\forall z\in F\forall y\in Y(yR'_\leq xz\lra yR'xz)$\\
\hskip12mm\=iff\hskip3mm\=  $\forall x\in A\forall z\in F(R'_\leq xz\subseteq R'xz)$\\
\>iff\> $\forall x\in A\forall z\in F(R^{111}xz\subseteq R''_\leq xz)$\\
\>iff\> $\forall x\in A\forall z\in F(R^{111}xz\subseteq R_\leq xz)$
\end{tabbing}
Hence, the hypothesis $R^{111}xz\subseteq R_\leq xz$, for all $x,z\in X$, implies $A\bigovert F\subseteq A\cap F$.

For the converse inclusion we use the assumption that for any $x,z\in X$ we have $R_\leq xz\subseteq R^{111}xz$, we let $u\in A\cap F$ and we show that $u\in A\bigovert F$.

Since $A=\bigcup_{x\in A}\Gamma x$ and $F=\bigcup_{z\in F}\Gamma z$, we obtain that $u\in A\cap F$ iff there exist elements $x\in A, z\in F$ such that $u\in \Gamma x$ and $u\in\Gamma z$. But $u\in\Gamma x\cap\Gamma z$ implies, given our case hypothesis and part (1) of this Proposition, that $u\in\Gamma x\bigovert\Gamma z$ and since $A\bigovert F=\bigvee_{x\in A}^{z\in F}\Gamma x\bigovert\Gamma z$ it follows then that $u\in A\bigovert F$. Therefore we also get $A\cap F\subseteq A\bigovert F$ and hence it follows that $A\cap F=A\bigovert F$.
\end{proof}

\begin{defn}\rm
\label{Heyting frame defn}
A {\em Heyting frame} is
an implicative frame $\mathfrak{F}=(X,\upv,Y,T^{\partial 1\partial})$ where $R^{111}xz= R_\leq xz$, for all $x,z\in X$.
\end{defn}
Thus, by Proposition \ref{heyting prop}, Heyting frames are exactly the implicative frames $\mathfrak{F}=(X,\upv,Y,T^{\partial 1\partial})$ whose full complex algebra $\mathfrak{F}^+$ is a (complete) Heyting algebra in which the stable set operation $\Ra$ induced by the frame relation $T^{\partial 1\partial}$ coincides with the residual of intersection in $\gpsi$.

We conclude this section by showing that if the frame is a Heyting frame, then the implication operation in its dual full complex algebra interprets intuitionistic implication in the received way.

\begin{prop}\rm
Let $\mathfrak{F}=(X,\upv,Y,T^{\partial 1\partial})$ be a Heyting frame. Then for any Galois stable sets $A,C$ and element $x\in X$
\[
\mbox{$x\in (A\Ra C)$ iff $\forall z\in X(x\leq z\lra(z\in A\lra z\in C))$.}
\]
\end{prop}
\begin{proof}
Let $\alpha_R$ be the image operator on subsets of $X$ generated by a ternary relation $R$, $\alpha_R(U,W)=\bigcup_{x\in U, z\in W}Rxz$.
Since $\alpha_R$ distributes over arbitrary unions, it is residuated in the powerset algebra with a map $\beta_R$, which is then defined  by $\beta_R(U,V)=\bigcup\{W\subseteq X\midsp \alpha_R(U,W)\subseteq V\}$.

It was shown in \cite[Theorem~3.14]{duality2}, that if the closure $\overline{\alpha}_R$ of the restriction of $\alpha_R$ to Galois sets is residuated, then its residual is the restriction $\beta_{R/}$ of $\beta_R$ on Galois sets, explicitly defined by $\beta_{R/}(A,C)=\bigcup\{F\in\gpsi\midsp\alpha_R(A,F)\subseteq C\}$. We have further shown in \cite[Lemma~3.15]{duality2} that $\beta_{R/}$ is equivalently defined by $\beta_{R/}(A,C)=\{u\in X\midsp \alpha_R(A,\Gamma u)\subseteq C\}$.

In our case of interest, $R=R^{111}=R_\leq$, $\overline{\alpha}_R(A,C)=\alpha_R(A,C)=A\cap C$ and, since we assume the frame is a Heyting frame, $\bigovert=\cap$ is indeed residuated with the implication operation $\Ra$ defined using the frame relation $T^{\partial 1\partial}$ (Definition~\ref{implication def}). By \cite[Lemma~3.15]{duality2} applied to our case of interest we obtain that $A\Ra C=\{x\in X\midsp A\cap \Gamma x\subseteq C\}$. But this means exactly that $x\in (A\Ra C)$ iff for any $z\in X$, if $x\leq z$ and $z\in A$, then $z\in C$, which is equivalent to the membership condition in the statement of the Proposition.
\end{proof}

\section{Choice-free Topological Duality}
\label{rep section}
\subsection{Integral Implicative Lattices}
\label{rep integral section}
The lattice representation we present is that of \cite{sdl}, recast in a choice-free manner in \cite[Theorem~4.6]{choiceFreeStLog} by switching from a Stone to a spectral topology.

\begin{thm}[Choice-free Lattice Representation]\label{choice-free lat rep}
  \rm
  Let $\mathbf{L}=(L,\leq,\wedge,\vee,0,1)$ be a bounded lattice and $(X,\upv,Y)$ its dual filter-ideal frame ($X=\filt(\mathbf{L}),Y=\idl(\mathbf{L})$), with ${\upv}\subseteq X\times Y$ defined by $x\upv y$ iff $x\cap y\neq\emptyset$. Let $\mathfrak{X}=(X,\mathcal{B})$ and $\mathfrak{Y}=(Y,\mathcal{C})$ be the spectral spaces generated by the bases $\mathcal{B}=\{X_a\midsp a\in L\}$ and $\mathcal{C}=\{Y^a\midsp a\in L\}$, respectively.

  Then the map $a\mapsto X_a$ is a lattice isomorphism $\mathbf{L}\iso{\tt KO}\mathcal{G}(\filt(\mathbf{L}))$ and the map $a\mapsto Y^a$ is a dual isomorphism $\mathbf{L}^\partial\iso{\tt KO}\mathcal{G}(\idl(\mathbf{L}))$.\telos
\end{thm}
We sketch the proof below, referring to \cite{duality2,choiceFreeStLog} for details.

As detailed in \cite{choiceFreeStLog}, the topological spaces $\mathfrak{X}=(X,\mathcal{B})$ and $\mathfrak{Y}=(Y,\mathcal{C})$ are spectral spaces \cite[Proposition~4.3]{choiceFreeStLog}, each of $\mathcal{B,C}$ is a meet semilattice ($X_a\cap X_b=X_{a\wedge b}$ and $Y^a\cap Y^b=Y^{a\vee b}$), each consists of the compact-open Galois stable and co-stable, respecitvely, sets and the Galois connection induced by $\upv$ restricts to a dual isomorphism $\mathcal{B}\iso\mathcal{C}^\mathrm{op}$, hence they are both (dually isomorphic) lattices \cite[Proposition~4.5]{choiceFreeStLog}.

Furthermore, the map $a\mapsto X_a$ is a lattice isomorphism $\mathbf{L}\iso{\tt KO}\mathcal{G}(\filt(\mathbf{L}))$ and the map $a\mapsto Y^a$ is a dual isomorphism $\mathbf{L}^\partial\iso{\tt KO}\mathcal{G}(\idl(\mathbf{L}))$ \cite[Proposition~4.6]{choiceFreeStLog}. The representation of \cite{choiceFreeStLog} differs from that we gave in \cite{sdl} only by resorting to a spectral, rather than a Stone topology and by restricting to proper filters/ideals only. The latter choice does not affect the argument and by the proof in \cite[Proposition~2.6]{mai-grishin} the representation of \cite{sdl} (hence that of \cite{choiceFreeStLog} that we use here) is a canonical extension of the lattice.

To represent implication, given an implicative lattice $\mathbf{L}=(L,\leq,\wedge,\vee,0,1,\ra)$, we apply the framework of \cite{duality2} where for each normal lattice operator a relation is added to the frame, such that if $\delta=(i_1,\ldots,i_n;i_{n+1})$ is the distribution type of the operator, then $\sigma=(i_{n+1};i_1\cdots i_n)$ is the sort type of the relation. In our case of interest in the present article and given that $\delta(\ra)=(1,\partial;\partial)$, we add a ternary relation $T^{\partial 1\partial}$ of the indicated sort type.

The canonical ternary relation $T=T^{\partial 1\partial}\subseteq Y\times (X\times Y)$ is defined as in \cite[Section~4.1]{duality2} using a point operator ${\leadsto}:X\times Y\lra Y$, where for $x\in X$ and $y\in Y$ we define $x\leadsto y=\bigvee\{y_{a\ra b}\midsp a\in x, b\in y\}$. We make the convention to designate principal filters by $x_a=a{\uparrow}$ and principal ideals by $y_a=a{\downarrow}$, so that $y_{a\ra b}$ is the principal ideal generated by an implication element $a\ra b$. Observe that if $a\in x\cap y$, then (since we assume the lattice to be integral) $(a\ra a)=1\in(x\leadsto y)$, which is an ideal, so that we obtain $(x\leadsto y)=\omega$ is the improper ideal. We let also $T^{11\partial}=T'$ be the Galois dual relation of $T$. We display the definitions below, for ease of reference, together with equivalent definitions for each of $T, T'$, using \cite[Lemma~4.4, Lemma~4.5]{duality2}.
\begin{tabbing}
\hskip8mm\=$x\leadsto v$\hskip8mm\==\hskip5mm\= $\bigvee\{y_{a\ra b}\midsp a\in x\in X\mbox{ and }b\in v\in Y\}\in Y=\idl(\mathbf{L})$\\
\>$yT^{\partial 1\partial}xv$\>iff\> $(x\leadsto v)\subseteq y$\hskip5mm\= iff\hskip4mm\= $\forall a,b(a\in x\;\wedge\;b\in v\lra(a\ra b)\in y)$\\
\>$uT'xv$\>iff\> $u\upv(x\leadsto v)$\> iff\> $\exists a,b(a\in x\;\wedge\;b\in v\;\wedge\;(a\ra b)\in u)$
\end{tabbing}

\begin{prop}\rm
\label{reduct prop}
The canonical frame $\mathbf{L}_+=\mathfrak{F}=(\filt(\mathbf{L}),\upv,\idl(\mathbf{L}),T^{\partial 1\partial})$ of an implicative lattice $\mathbf{L}=(L,\leq,\wedge,\vee,0,1,\ra)$ is an implicative frame (in the sense of Definition \ref{implicative frame defn}). Moreover, if the lattice is distributive, then intersection distributes over arbitrary joins in $\gpsi$, hence the full complex algebra of the frame is a (complete) residuated Heyting algebra.
\end{prop}
\begin{proof}
Axioms (F1)--(F4) were verified more generally for any normal lattice expansion in \cite[Lemma~4.3, Lemma~4.6]{duality2}. For axiom (F0), if $x$ is a filter, $y$ an ideal and $x\upv y$, let $a\in x\cap y$. Then $a\ra a=1\in(x\leadsto y)$, so that $(x\leadsto y)=\omega$ is the improper ideal (the whole lattice) and then for any filter $u$ we have $u\upv(x\leadsto y)$, i.e. $uT^{11\partial}xy$.

Conversely, suppose that for all filters $u$ we have $u\upv(x\leadsto y)$. Then considering the trivial filter $u=\{1\}$ we conclude that $1\in(x\leadsto y)$. By definition of $x\leadsto y$, there exist $a\in x, b\in y$ such that $1\leq a\ra b$. By the integrality axiom in the lattice, this is equivalent to $a\leq b$. Then $a\in x\cap y$, i.e. $x\upv y$.

The case of a distributive lattice was presented in \cite[Theorem~5.2]{choiceFreeStLog}, but we repeat it here as this article has not at this point appeared in print.

To prove that $\gpsi$ is a completely distributive lattice, note first that both lattice join $\vee$ and meet $\wedge$ are trivially normal lattice operators in the sense of Definition \ref{normal lattice op defn}, but meet is an operator (in the J\'{o}nsson-Tarski sense) only when it distributes over joins. When this is the case, meet also has the distribution type $(1,1;1)$. Its $\sigma$-extension $\wedge_\sigma$, is constructed as outlined in \cite[Section~4.1]{duality2}. Specifically, letting $\wedge=f$, the point operator $\widehat{f}$ on filters is defined by $\widehat{f}(x,z)=\bigvee\{x_{a\wedge b}\midsp a\in x\mbox{ and }b\in z\}$ and the canonical relation $R_\wedge$ is then defined by $xR_\wedge uz$ iff $\forall a,b(a\in u\mbox{ and } b\in z\lra a\wedge b\in x)$, using Lemma~4.4 of \cite{duality2}.  Note that $R_\wedge$ is the upper bound relation of Proposition \ref{upper bound rel prop}. Considering the image operator $\alpha_R:\powerset(X)\times\powerset(X)\lra\powerset(X)$ defined by $\alpha_R(U,W)=\{x\in X\midsp \exists u\in U\exists z\in W\; xR_\wedge uz\}$,  we obtain that  $\alpha_R(A,C)=A\cap C$, for $A,C\in \gpsi$. By \cite[Lemma~4.6]{duality2}, all sections of the Galois dual relation of $R_\wedge$ are stable. It then follows by Proposition \ref{upper bound rel prop} that intersection distributes over arbitrary joins, in other words, $\gpsi$ is a completely distributive lattice. Combining with Proposition~\ref{residuation in frame prop} it is concluded that the full complex algebra of the frame is a residuated Heyting algebra.
\end{proof}

\begin{prop}\rm
The lattice representation map $a\mapsto X_a$ is an isomorphism $\mathbf{L}\iso\type{KO}\mathcal{G}(X)$ of integral implicative lattices.
\end{prop}
\begin{proof}
The image $X_a$ of a lattice element $a$ is the set of all filters containing the principal filter $x_a$, hence $X_a=\Gamma x_a$, which is a clopen element because $\Gamma x_a={}\rperp\{y_a\}$, where $y_a$ is the principal ideal generated by $a$. This is because $x\upv y_a$ iff $x\cap y_a\neq\emptyset$ iff $a\in x$. Therefore $X_a\Ra X_b=\Gamma x_a\Ra{}\rperp\{y_b\}=\{u\in X\midsp uT'x_ay_b\}$, by Proposition \ref{Ra prop}. In the canonical frame $uT'x_ay_b$ holds iff $u\upv(x_a\leadsto y_b)$ iff $\exists e,d(e\in x_a,d\in y_b$ and $(e\ra d)\in u)$, using \cite[Lemma~4.5]{duality2}. Given the monotonicity properties of implication, this means that $uT'x_ay_b$ iff $(a\ra b)\in u$. Conclude from this that $X_a\Ra X_b=X_{a\ra b}$ and this proves the claim that the representation map is a homomorphism of implicative lattices. Bijectivity was shown in the latice representation result  \cite[Theorem~5.5]{duality2} (or its choice-free version \cite[Theorem~4.6]{choiceFreeStLog}).
\end{proof}

We conclude this section by providing a proof of Theorem~\ref{reduct thm}, stated below as a corollary to the results already obtained.

\begin{coro}\rm\label{reduct coro}
Every integral implicative lattice $\mathbf{L}$ is a reduct of an integral residuated lattice $\mathbf{R}$. If $\mathbf{D}$ is a distributive integral implicative lattice, then it is a reduct of a residuated Heyting algebra $\mathbf{H}$. In both cases, $\mathbf{R,H}$ can be constructed as a canonical extension of $\mathbf{L}$ and of $\mathbf{D}$, respectively.
\end{coro}
\begin{proof}
The claim is proven in Proposition~\ref{reduct prop}. That the representation is a canonical extension of the represented implicative lattice follows by the fact that the lattice representation delivers a canonical extension, by \cite[Proposition~2.6]{mai-harding}. In \cite[Proposition~28]{kata2z} we have detailed the proof that our canonical representation of a normal lattice operator delivers its $\sigma$-extension (in the terminology of \cite[Section~4]{mai-harding}) if its output type is 1, and it delivers its $\pi$-extension when the output type is $\partial$.
\end{proof}

\subsection{Heyting Algebras}
\label{rep Heyting section}
Assume $\mathbf{H}=(H,\leq,\wedge,\vee,0,1,\ra)$ is a Heyting algebra. By Proposition~\ref{reduct prop}, the full complex algebra $\mathfrak{F}^+$ of its dual frame $\mathfrak{F}=(\filt(\mathbf{H}),\upv,\idl(\mathbf{H}),T)$, where $T=T^{\partial 1\partial}\subseteq Y\times(X\times Y)$ and $X=\filt(\mathbf{H}),Y=\idl(\mathbf{H})$, is a residuated Heyting algebra with $\Ra$ the residual of $\cap$, defined by
\begin{equation}
\label{implication residual equation}
A\Ra C=\{u\in X\midsp\forall z\in X(z\in A\;\wedge\;u\leq z\lra z\in C)\},
\end{equation}
and $\bigovert_T,\Ra_T$ the residuated operators induced by the frame relation $T$, as shown in Proposition~\ref{Ra prop}, defined by
\begin{equation}\label{T-implication residual equation}
A\Ra_TC=\{u\in X\midsp \forall x\in X\forall y\in Y(x\in A\;\wedge\; C\upv y\lra uT'xy)\}.
\end{equation}

\begin{prop}\rm
\label{heyting frame prop}
If $\mathbf{H}=(H,\leq,\wedge,\vee,0,1,\ra)$ is a Heyting algebra, then its canonical frame $\mathfrak{F}=(\filt(\mathbf{H}),\upv,\idl(\mathbf{H}),T)$, where $T=T^{\partial 1\partial}\subseteq Y\times(X\times Y)$ and $X=\filt(\mathbf{H}),Y=\idl(\mathbf{H})$, is a Heyting frame (Definition \ref{Heyting frame defn}).
\end{prop}
\begin{proof}
To prove that the canonical frame of $\mathbf{H}$ is a Heyting frame it suffices to show  that $\Ra\; =\; \sub{\Ra}{T}$, or that $\cap=\sub{\bigovert}{T}$, or that $R_\wedge xz = R^{111}xz$ for all $x,z\in X$, given Proposition \ref{heyting prop} and uniqueness of adjoints. Furthermore, given the distribution properties of $\Ra$ and $\sub{\Ra}{T}$ it suffices to prove that for any $x\in X$ and $y\in Y$ we have $\Gamma x\Ra{}\rperp\{y\}=\Gamma x\sub{\Ra}{T}{}\rperp\{y\}$. Equivalently, by part 3 of Proposition~\ref{Ra prop}, it suffices to show that $uT^{11\partial}xy$ iff $u\in\Gamma x\Ra{}\rperp\{y\}$. Using equation \eqref{implication residual equation} and the fact that $T^{11\partial}$ is the Galois dual relation of $T^{\partial 1\partial}$ and that in the canonical frame $uT'xy$ is equivalent to $u\upv(x\leadsto y)$, the claim further reduces to showing that $\Gamma x\cap\Gamma u\subseteq{}\rperp\{y\}$ iff $u\upv (x\leadsto y)$. Furthermore, in the lattice of filters $\Gamma u\cap\Gamma x=\Gamma(x\vee u)$ so that $\Gamma u\cap\Gamma x\subseteq{}\rperp\{y\}$ iff $(x\vee u)\upv y$.

By the above, it suffices to prove that if $x,u\in\filt(\mathbf{H})$ and $y\in\idl(\mathbf{H})$, then $u\upv (x\leadsto y)$ iff $(x\vee u)\upv y$.

Assume $u\upv (x\leadsto y)$ and let then $a\in x, b\in y$ such that $(a\ra b)\in u$. Both $a,a\ra b$ are in the join $x\vee u$ of the filters $x,u$, hence by $a\wedge(a\ra b)\leq b$ we also get $b\in (x\vee u)$. Then $b\in (x\vee u)\cap y\neq\emptyset$ and so by definition $(x\vee u)\upv y$.

Conversely, assume $(x\vee u)\upv y$ and let $e\in (x\vee u)\cap y$. Let $a\in x, b\in u$ such that $a\wedge b\leq e$.  By residuation in $\mathbf{H}$, $b\leq a\ra e$. Since $a\in x$ and $e\in y$, we have $(a\ra e)\in (x\leadsto y)$. Since $b\in u$ and $b\leq a\ra e$ we also have $(a\ra e)\in u$. Hence $u\cap(x\leadsto y)\neq\emptyset$ which means, by definition, $u\upv (x\leadsto y)$.
\end{proof}
The following has then been established.
\begin{coro}[Choice-free Representation of Heyting Algebras]\rm
\label{Heyting rep coro}
If $\mathbf{H}$ is a Heyting algebra, then the full complex algebra $\mathfrak{F}^+$ (the algebra $\gpsi$ of stable sets of filters) of its dual frame $\mathfrak{F}=(\filt(\mathbf{H}),\upv,\idl(\mathbf{H}),T)$ is a complete Heyting algebra and a canonical extension of $\mathbf{H}$, which is identified as the subalgebra $\type{KO}\mathcal{G}(X)$ of compact-open Galois stable sets. Implication, the residual of intersection, is defined by equation \eqref{implication residual equation}, which is equivalent to the definition by equation \eqref{T-implication residual equation}.\telos
\end{coro}

\subsection{Duality}
\label{duality section}
The (choice-free) representation of integral implicative lattices (and of Heyting algebras, in particular) extends to a full functorial duality. We sketch the argument rather than discussing it in any detail, as there is really nothing new to add, except for specializing to the case of interest in this paper a result that has been thoroughly presented in \cite{duality2}, also in \cite{choiceFreeStLog}, to which we refer the interested reader for details.

In Section \ref{impl lattice section} we defined the category $\mathbf{NLE}_{1\partial\partial}$ of integral implicative lattices. In Section \ref{complex algebra section} and, in particular, in Definition \ref{implicative frame defn} we defined implicative frames, the objects of the category $\mathbf{SRF}_{1\partial\partial}$ of sorted residuated frames $\mathfrak{F}=(X,\upv,Y,T)$, with a ternary relation $T=T^{\partial 1\partial}\subseteq Y\times(X\times Y)$. Morphisms of the frame category are the weak bounded morphisms of \cite[Definition~3.20]{duality2}. In \cite[Section~4]{duality2}  a contravariant functor $\type{F}:\mathbf{NLE}_\tau\lra\mathbf{SRF}_\tau^\mathrm{op}$, was defined by a canonical frame construction, which specializes to our present case of implicative lattices and frames, verifying in \cite[Proposition~4.9]{duality2} that the duals of lattice expansion homomorphisms are weak bounded morphisms. For duality purposes, a smaller category $\mathbf{SRF}^*_\tau$ was defined, axiomatized in \cite[Table~3]{duality2}. This was modified in \cite[Table~2]{choiceFreeStLog}, switching from a Stone to a spectral topology in order to obtain a choice-free result. A contravariant functor $\type{L}^*:\mathbf{SRF}_\tau^*\lra\mathbf{NLE}_\tau^\mathrm{op}$ was defined in \cite{duality2}, specializing in our present case to a functor $\type{L}^*:\mathbf{SRF}_{1\partial\partial}^*\lra\mathbf{NLE}_{1\partial\partial}^\mathrm{op}$, mapping a frame to the subalgebra of compact-open Galois stable sets of its full complex algebra \cite[Proposition~5.2]{duality2} and a weak bounded morphism to a homomorphism of the dual lattice algebras \cite[Proposition~5.3]{duality2}. Theorem 5.8 of \cite{duality2} concluded with the proof of duality, which specializes to (choice-free) dualities for implicative lattices and frames, as well as for Heyting algebras and their dual Heyting frames.

\section{Applications and Further Issues}
\label{conc section}
The received relational semantics for $\mathbf{IPC}$ uses frames $(W,R)$, where $R\subseteq W\times W$ is a reflexive and transitive relation, sentences are interpreted as $R$-closed sets, with implication interpreted by the clause
\begin{equation}\label{int clause}
w\forces\varphi\ra\psi\;\mbox{ iff }\; \forall u(wRu\;\mbox{ $\lra$}\;(u\forces\varphi\lra u\forces \psi))
\end{equation}
The box operator generated by $R$, $[R]U=\{w\in W\midsp\forall u(wRu\lra u\in U)\}$, for $U\subseteq W$, is an $\mathbf{S4}$-modality (a $\mathbf{KT4}$-modality) under the given assumptions for the accessibility relation $R$.
Subintuitionistic logics \cite{Visser1981,dosen-modal-trans-in-KD,corsi-subint,restall-subint,Wansing1997-subint,celani-subint-2001,weak-subint} arise by tampering with the properties of the accessibility relation $R$, as originally proposed by Corsi, Restall and Wansing, hence modifying the logic of $[R]$ to $\mathbf{KT, KD,KTB}$ etc (see \cite[Section~6, Lemma~6.1]{corsi-subint}) or, equivalently, by considering logics in the language of $\mathbf{IPC}$ whose modal companions under the GMT translation \cite{goedel} are weaker than $\mathbf{S4}$, a direction explicitly taken by Do\v{s}en.

Dropping residuation of implication with conjunction is the distinctive mark of logics weaker than $\mathbf{IPC}$. This is not related to subintuitionistic systems alone, but it includes all substructural logics, as well. Unlike subintuitionistic logics, which are defined by weakening frame conditions in intuitionistic frames, substructural logics have been defined proof-theoretically, by weakening the proof system of $\mathbf{IPC}$, removing one or more of the so-called structural rules of association, exchange, weakening and contraction. There has been no attempt, to the best of this author's knowledge, to define substructural logics by the same process of weakening frame conditions, as for subintuitionistic logics. Evidently, in such an approach, the set-theoretic (relational) semantics of $\mathbf{IPC}$ is reconsidered, and we proposed sorted frame semantics for this purpose. We carried this project out in this article and in \cite{redm}, but without tampering with the semantics of intuitionistic implication. In particular, though we have not considered this yet in due detail, we can regard distributive `subintermediate' logics (subintuitionistic, or substructural) as  fragments of the logic of residuated Heyting algebras, which are like the residuated Boolean algebras of J\'{o}nsson and Tsinakis\cite{residBA}, except for weakening the underlying Boolean to a Heyting algebra structure.

Sorted frame semantics for $\mathbf{IPC}$ can be extended to include a treatment of Intuitionistic Modal Logic $\mathbf{IML}$, combining \cite{pnsds,duality2} with the present article, but we leave this for future research.


\end{document}